\documentclass[journal]{IEEEtran}
\usepackage{amsmath,amssymb,amsthm,bbm,color,abbrevs,subfigure,graphicx,algorithm,algpseudocode,url}
\usepackage[usenames,dvipsnames]{xcolor}
\usepackage[normalem]{ulem}

\newcommand \E {\mathbbm{E}}
\newcommand \diag{\mathrm{dg}}
\newcommand \argmin{\mathrm{argmin}}

\newcommand \bp {\mathbf{p}}
\newcommand \bq {\mathbf{q}}
\newcommand \bv{\mathbf{v}}
\newcommand \bxi{\boldsymbol{\xi}}
\newcommand \bz{\mathbf{z}}
\newcommand \by{\mathbf{y}}
\newcommand \blambda{\boldsymbol{\lambda}}
\newcommand \bnu{\boldsymbol{\nu}}
\newcommand \bmu{\boldsymbol{\mu}}
\newcommand \bepsilon{\boldsymbol{\epsilon}}
\allowdisplaybreaks

\begin{document}
\title{Two-Timescale Stochastic Dispatch\\
of Smart Distribution Grids}
\author{
	Luis M. Lopez-Ramos,~\IEEEmembership{Student Member,~IEEE,}
	Vassilis Kekatos,~\IEEEmembership{Member,~IEEE,}\\
	Antonio G. Marques,~\IEEEmembership{Senior Member,~IEEE,} and
	Georgios B. Giannakis,~\IEEEmembership{Fellow,~IEEE}
	

\thanks{This work was supported by the Spanish Ministry of Education FPU Grant AP2010-1050; CAM Grant S2013/ICE-2933; MINECO Grant TEC2013-41604-R; and NSF grants 1423316, 1442686, 1508993, and 1509040.}

\thanks{L.-M. Lopez-Ramos and Antonio G. Marques are with the Dept. of Signal Theory and Communications, King Juan Carlos Univ., Fuenlabrada, Madrid 28943, Spain. V. Kekatos is with the ECE Dept., Virginia Tech, Blacksburg, VA 24061, USA. G. B. Giannakis is with the Digital Technology Center and the ECE Dept., University of Minnesota, Minneapolis, MN 55455, USA. Emails: luismiguel.lopez@urjc.es, kekatos@vt.edu, antonio.garcia.marques@urjc.es, georgios@umn.edu.}

}
\markboth{Lopez-Ramos, Kekatos, Marques, and Giannakis: Two-Timescale Stochastic Dispatch of Smart Distribution Grids}{}

\maketitle

\begin{abstract}
Smart distribution grids should efficiently integrate stochastic renewable resources while effecting voltage regulation. The design of energy management schemes is challenging, one of the reasons being that energy management is a multistage problem where decisions are not all made at the same timescale and must account for the variability during real-time operation. The joint dispatch of slow- and fast-timescale controls in a smart distribution grid is considered here. The substation voltage, the energy exchanged with a main grid, and the generation schedules for small diesel generators have to be decided on a slow timescale; whereas optimal photovoltaic inverter setpoints are found on a more frequent basis. While inverter and looser voltage regulation limits are imposed at all times, tighter bus voltage constraints are enforced on the average or in probability, thus enabling more efficient renewable integration. Upon reformulating the two-stage grid dispatch as a stochastic convex-concave problem, two distribution-free schemes are put forth. An average dispatch algorithm converges provably to the optimal two-stage decisions via a sequence of convex quadratic programs. Its non-convex probabilistic alternative entails solving two slightly different convex problems and is numerically shown to converge. Numerical tests on a real-world distribution feeder verify that both novel data-driven schemes yield lower costs over competing alternatives.
\end{abstract}


\section{Introduction}\label{s:intro}
With increasing renewable generation, energy management of power distribution grids is becoming a computationally challenging task. Solar energy from photovoltaic (PV) units can change significantly over one-minute intervals. The power inverters found in PV units can be commanded to curtail active power generation or adjust their power factor within seconds~\cite{nrel2008liub}, \cite{Carvalho}. At a slower timescale, distribution grid operators exchange energy with the main grid hourly or on a 10-minute basis, and may experience cost penalties upon deviating from energy market schedules~\cite{VaWuBi11}. Moreover, voltage regulation equipment and small diesel generators potentially installed in microgrids respond at the same slower timescale. As a result, comprehensive designs to optimize such diverse tasks call for \emph{multistage} smart grid dispatch solutions.

Spurred by demand-response programs and the use of PV inverters to accomplish various grid tasks~\cite{Turitsyn11}, \emph{single-stage} dispatch schemes for distribution grids have been an active area of research. Power inverters can be controlled using localized rules for voltage regulation, see e.g.,~\cite{ZDGT13,VKZG15,Bolognani15,Zhu16}. Assuming two-way communication between buses and the utility operator, dispatching a distribution system can be posed as an optimal power flow (OPF) problem. Centralized schemes use nonlinear program solvers~\cite{Canizares}; or rely on convex relaxations of the full ac model of balanced~\cite{FL1,tac2015gan}, or unbalanced grids~\cite{tse2014emiliano}. Distributed solvers with reduced computational complexity have been devised in~\cite{EmilSG13,Baosen14,PengLow16}.

Nevertheless, the efficient and secure operation of distribution grids involves decisions at different timescales. A dynamic programming approach for a two-stage dispatch is suggested in~\cite{FCL}: The taps of voltage regulators are set on a slow timescale and remain fixed for consecutive shorter time slots over which elastic loads are dispatched; yet the flexibility of loads is assumed known \textit{a priori}. Alternatively, centrally computed OPF decisions can be communicated to buses at a slow timescale, while on a faster timescale, PV power electronics are adjusted to optimally track variations in renewable generation and demand~\cite{hierarchies,Emil16}. 
Relying on approximate grid models, the latter schemes yield a fully localized real-time implementation. However, they presume smooth system transitions and dispatch slow-responding units for a single deterministic fast-timescale scenario. 

Multistage dispatching under uncertainty is routinely used in transmission systems and microgrids~\cite{conejo2010decision}. Robust approaches find optimal slow-timescale decisions for the worst-case fast-timescale outcome; see~\cite{YGG13} and references therein. To avoid the conservativeness of robust schemes, probabilistic approaches postulate a probability density function (pdf) for demand, wind generation, and system contingencies to find day-ahead grid schedules~\cite{bouffard2005market}, \cite{Bienstock}. The risk-limiting dispatch framework adjusts multistage decisions as the variance of the random variables involved decreases while approaching actual time~\cite{VaWuBi11}. Decisions can be efficiently calculated 
only for convenient pdfs
for a network-constrained risk-limiting dispatch under additional 
transmission congestion assumptions~\cite{rajagopal2014network}. As a third alternative, stochastic sample approximation approaches yield optimal slow-timescale decisions using samples drawn from the postulated pdf; see e.g.,~\cite{BoGa08}, \cite{ZhangGreenTech}. 

Returning to distribution grids, PV inverters could be overloaded sporadically in time and across buses to accommodate solar fluctuations and prevent overvoltages~\cite{gtd2013kirtley}. The spatiotemporal overloading of power system components (such as inverters, bus voltages, line flows) could thus constitute an additional means for integrating renewables in smart grids. Nonetheless, ensuring that overloading occurs sparingly couples decisions across time. The single-stage scheme of \cite{ergodic} finds optimal PV setpoints while limiting time averages of overloaded quantities. The latter approach has been also adopted in \cite{G-SIP15} for dispatching a transmission system in a day-ahead/real-time market setup under load shedding.

Jointly dispatching slow- and fast-timescale distribution grid resources under \emph{average} or \emph{probabilistic} constraints over fast-timescale decisions is considered here. Our contributions are three-fold. First, using an approximate grid model, the expected cost over a slow control period is minimized while inverter and looser bus voltage constraints are satisfied at all fast-timescale slots. Further, tighter voltage limits are enforced either on the average or in probability across successive fast-timescale slots (see Sec.~\ref{s:problem}). Two-stage grid dispatch is formulated as a convex-concave optimization in Sec.~\ref{s:analysis}. Second, adopting a stochastic saddle-point approximation scheme from \cite{Nemirovski09}, the provably convergent algorithm of Sec.~\ref{s:ada} finds the optimal slow-timescale decisions in the case of average constraints. Third, for the case of non-convex \textit{probabilistic} constraints, an algorithm solving two similar convex problems for each fast-timescale period is put forth in Sec.~\ref{s:pda}. Albeit the related expected recourse function enjoys zero-duality gap~\cite{ribeiro10}, the overall two-stage dispatch is not convex-concave; and the algorithm's performance is only numerically validated. Both algorithms require only samples rather than pdfs of loads and solar generation, and involve solving simple convex quadratic programs. Numerical tests in Sec.~\ref{s:numerical} on a 56-bus feeder corroborate the validity of our findings.

Regarding \emph{notation}, lower-(upper-)case boldface letters denote column vectors (matrices), with the only exception of the power flow vectors, which are uppercase. Calligraphic letters are used to denote sets. Symbol $^\top$ denotes transposition, while $\mathbf{0}$ and $\mathbf{1}$ are the all-zeros and all-ones vectors of appropriate dimensions. The indicator function $\mathbbm{1}{\{\cdot\}}$ equals 1 when its argument is true, and 0 otherwise. A diagonal matrix with the entries of vector $\mathbf{x}$ on its main diagonal is denoted by $\diag(\mathbf{x})$. The operator $[\cdot]_+$ projects its argument onto the positive orthant; $\E [\cdot]$ denotes expectation and $\Pr\{\cdot\}$ probability.

\section{Problem Formulation}\label{s:problem}
Consider a distribution grid whose energy needs are procured by distributed renewable generation, distributed conventional (small diesel) generators, and the main grid. The distribution grid operator aims at serving load at the minimum cost while respecting voltage regulation and network constraints. Energy is exchanged with the main grid at whole-sale electricity prices through the feeder bus. To effectively integrate stochastic renewable generation, the focus here is on short-term grid dispatch. To that end, the distribution grid is operated at two timescales: a slower timescale corresponds to 5- or 10-min real-time energy market intervals, while the inverters found in PVs are controlled at a faster timescale of say 10-sec intervals. One period of the slower timescale is comprised by $T$ faster time slots indexed by $t=1,\ldots,T$.

The grid is operated as a radial network with $N+1$ buses rooted at the substation bus indexed by $n=0$. The distribution line feeding bus $n$ is also indexed by $n$ for $n=1,\ldots,N$. Let $p_{n,t}$ and $q_{n,t}$ denote respectively the net active and reactive power injections at bus $n$ and slot $t$; the $N$-dimensional vectors $\bp_t$ and $\bq_t$ collect the net injections at all buses except for the substation. Diesel generators are dispatched at the slower timescale to generate $\bp^{d}$ throughout the subsequent $T$ slots at unit power factor. During slot $t$, PVs can contribute solar generation up to $\overline{\mathbf{p}}_t^r$ that is modeled as a random process. Smart inverters perform active power curtailment and reactive power compensation by following the setpoints $\bp^r_t$ and $\bq^r_t$ commanded by the utility operator. Load demands $\bp^l_t$ and $\bq^l_t$ are also modeled as random processes. To simplify the exposition, $(\bp^l_t,\bq^l_t)$ are assumed inelastic and known at the beginning of slot $t$; although elastic loads can be incorporated without any essential differences. The operator buys a power block $p_0^a$ from the main grid at the slow timescale, which can be adjusted to $p_{0,t}:=p_0^a + p^\delta_{0,t}$ in actual time.

Voltage regulation is effected by controlling (re)active power injections at slot $t$. Let $v_{n,t}$ denote the squared voltage magnitude at bus $n$ and slot $t$, and $\bv_t$ the vector collecting $\{v_{n,t}\}_{n=1}^N$. The substation voltage $v_0^a$ is controlled at the slower timescale \cite{FCL}, while voltage magnitudes at all buses must adhere to voltage regulation standards, e.g., ANSI C84.1 and EN50160 in \cite{ansic84}, \cite{en50160}. These standards differentiate between a narrower voltage regulation range denoted here by $\mathcal{V}_A$ in which voltages should lie most of the time; and a wider range $\mathcal{V}_B$ (with $\mathcal{V}_A\subset \mathcal{V}_B$) whom voltages should not exceed at any time. One of the goals of this work is to leverage this flexibility to design dispatch schemes that: i) guarantee that voltages lie in $\mathcal{V}_B$ at all times, while ii) they belong to $\mathcal{V}_A$ in a stochastic fashion. To this end, two alternative schemes are presented, the difference between them being how constraint ii) is formulated. The first scheme guarantees that the \textit{average} voltage lies in $\mathcal{V}_A$, whereas the second one maintains the \textit{probability} of under-/over-voltage at a specified low value.

\subsection{Grid modeling}
To capture voltage and network limitations, the distribution grid is captured by the approximate linear distribution flow (LDF) model, which is briefly reviewed next~\cite{BW3}. Let $\mathbf{r}$ and $\mathbf{x}$ be accordingly the vectors of line resistances and reactances across lines. Define also the branch-bus incidence matrix $\tilde{\mathbf{A}}\in\mathbb{R}^{N\times (N+1)}$ whose $(i,j)$-th entry is
\begin{equation} \label{eq:Amatrix}
\tilde{A}_{ij} = 
	\begin{cases}
		+1&,~\textrm{if $j-1$ is the source bus of line $i$}\\
		-1&,~\textrm{if $j-1$ is the destination bus of line $i$}\\
		0&,~\textrm{otherwise}.
	\end{cases}
\end{equation}
Partition $\tilde{\mathbf{A}}$ into its first column and the \emph{reduced} branch-bus incidence matrix $\mathbf{A}$ as $\tilde{\mathbf{A}}=[\mathbf{a}_0~\mathbf{A}]$. Ignoring line losses, the LDF model asserts that the vectors of active and reactive line power flows at time $t$ can be approximated by
\begin{equation}\label{eq:flows}
\mathbf{P}_t=\mathbf{F}^{\top}\bp_t~\textrm{and}~\mathbf{Q}_t=\mathbf{F}^{\top}\bq_t
\end{equation}
where $\mathbf{F}:=\mathbf{A}^{-1}$. Moreover, the squared voltage magnitudes can be expressed as~\cite{BW3}, \cite{FCL13}, \cite{SGC15}
\begin{equation}\label{eq:voltage}
\bv_t  = 2\mathbf{R}\bp_t + 2\mathbf{X}\bq_t + v_0^d \mathbf{1}
\end{equation}
where $\mathbf{R}:=\mathbf{F}\diag(\mathbf{r})\mathbf{F}^{\top}$ and $\mathbf{X}:= \mathbf{F}\diag(\mathbf{x})\mathbf{F}^{\top}$. Let us define the voltage regulation regions
\begin{subequations}\label{eq:vm}
\begin{align}
	\mathcal{V}_A&:=\left\{\mathbf{v}:~\underline{v}_{A}\mathbf{1}\leq \mathbf{v} \leq \overline{v}_{A}\mathbf{1}\right\}  \label{eq:vm:tight}\\
	\mathcal{V}_B&:=\left\{\mathbf{v}:~\underline{v}_{B}\mathbf{1} \leq \mathbf{v}\leq \overline{v}_{B}\mathbf{1}\right\} \label{eq:vm:loose}
\end{align}
\end{subequations}
with $\overline{v}_{B}\geq \overline{v}_{A}$ and $\underline{v}_{B}\leq \underline{v}_{A}$. Compliance with $\mathcal{V}_A$ can be imposed either on the average as $\E_t\left[\bv_t\right]\in \mathcal{V}_A$, or in probability as $\Pr\{\bv_t \in \mathcal{V}_A\} \geq 1- \alpha$ for some small $\alpha$. Either way, safe grid operation requires that $\mathbf{v}_t\in \mathcal{V}_B$ at all times $t$. Within the optimization horizon, the random processes involved (demand and renewable generation) can be assumed ergodic, i.e., their time averages converge to their ensemble averages. For this reason, voltage constraints pertaining to $\mathcal{V}_A$ will be referred to as \emph{ergodic}.

According to \eqref{eq:flows}, if $\mathbf{f}_n$ is the $n$-th column of $\mathbf{F}$, the squared power flow on line $n$ can be written as $P_{n,t}^2=\bp_t^{\top}\mathbf{f}_n\mathbf{f}_n^{\top}\bp_t$ and $Q_{n,t}^2=\bq_t^{\top}\mathbf{f}_n\mathbf{f}_n^{\top}\bq_t$. Imposing the upper limit $\overline{S}_n$ on the apparent flow on line $n$ is thus expressed as the convex quadratic constraint
\begin{equation}\label{eq:appline}
\bp_t^{\top}\mathbf{f}_n\mathbf{f}_n^{\top}\bp_t +  \bq_t^{\top}\mathbf{f}_n\mathbf{f}_n^{\top}\bq_t\leq \overline{S}_n^2.
\end{equation}
Assuming voltage magnitudes to be close to unity, active power losses can be approximated as $\sum_{n=1}^N r_n(P_{n,t}^2+Q_{n,t}^2)$~\cite{Chertkov}, which from \eqref{eq:flows}--\eqref{eq:voltage}, can be equivalently expressed as $\mathbf{P}_t^{\top}\diag(\mathbf{r})\mathbf{P}_t + \mathbf{Q}_t^{\top}\diag(\mathbf{r})\mathbf{Q}_t=\bp_t^{\top}{\mathbf{R}}\bp_t + \bq_t^{\top}{\mathbf{R}}\bq_t$. Thus, the active power injection at the substation is approximately
\begin{equation}\label{eq:activelosses}
p_{0,t} =-\mathbf{1}^\top \bp_t + \bp_t^\top{\mathbf{R}}\bp_t + \bq_t^\top{\mathbf{R}}\bq_t
\end{equation}

Regarding smart inverters, the tuple $(p_{n,t}^r, q_{n,t}^r)$, which denotes the power injection from the inverter located on bus $n$ at slot $t$, should belong to the feasible set
\begin{subequations} \label{def:omegaset}
\begin{align} 
\Omega_{n,t} :=\big\{(p_{n,t}^r, q_{n,t}^r):&~0\leq p_{n,t}^r\leq \overline{p}_{n,t}^r,\label{eq:omega:active}\\
&~|q_{n,t}^r|\leq \phi_n p_{n,t}^r,\label{eq:omega:power_factor}\\
&~(p_{n,t}^r)^2+ (q_{n,t}^r)^2\leq \overline{s}_n^2\big\}\label{eq:omega:apparent_power}
\end{align}
\end{subequations}
that is random and time-variant due to the variability of $\overline{p}^r_{n,t}$. 
Constraint \eqref{eq:omega:active} limits the active power generation according to the available solar power; constraint \eqref{eq:omega:power_factor} enforces the lower limit $\cos(\arctan (\phi_n))$ on the power factor (lagging or leading); and \eqref{eq:omega:apparent_power} limits the inverter apparent power.

\subsection{Operation costs} \label{ss:costs}
If PV owners are compensated at price $\boldsymbol{\pi}$ for the active power surplus they inject into the distribution grid, the related utility cost at slot $t$ is
$C_{\textrm{PV}}(\bp^r_t):= \boldsymbol{\pi}^\top [\bp^r_t - \bp^l_t]_+$ with $[\cdot]_{+}:=\max\{0,\cdot\}$ applied entrywise on vector $\bp^r_t - \bp^l_t$. The diesel generation cost is represented by $C_\textrm{D}(\bp^d)$. Regarding energy transactions with the main grid, the power block $p_0^a$ bought in advance is charged at a fixed and known price $\beta$. Deviating from $p_0^a$ by $p^\delta_{0,t}$ at slot $t$ is charged at
\begin{equation}\label{eq:RTcost}
C^t(p^\delta_{0,t}):= \gamma_b[p^\delta_{0,t}]_{+} - \gamma_s[-p^\delta_{0,t}]_{+}
\end{equation}
for known prices $(\gamma_b,\gamma_s)$. To avoid arbitrage, it is assumed that $0< \gamma_s < \beta < \gamma_b$; see e.g., \cite{VaWuBi11},~\cite{rajagopal2014network}. Then, the deviation charge can also be expressed as $C^t(p^\delta_{0,t})= \max\{\gamma_b p^\delta_{0,t}, \gamma_s p^\delta_{0,t} \}$, which is certainly convex~\cite{YGG13}.

\subsection{Optimal grid dispatch} \label{ss:statement}
Depending on the way compliance with voltage regulation region $\mathcal{V}_A$ is enforced, two grid dispatch formulations are developed next. Commencing with the \emph{average dispatch}, the optimal grid operation is posed as
\begin{subequations}\label{eq:p}
\begin{align}
\mathsf{P}_a^\ast := \min 	~&C_\textrm{D}(\bp^{d}) + \beta p_0^a +\E_t\left[C^t(p^\delta_{0,t}) + C_{\textrm{PV}}(\bp^r_t)\right] 					\label{eq:p:objective}			\\
\text{s.to:}~~& \bp_t = \bp^{r}_t - \bp^{l}_t + \bp^{d}									\label{eq:p:active_balance}		\\
				&  \bq_t =\bq^{r}_t - \bq^{l}_t 										\label{eq:p:reactive_balance}	\\	
				& p_{0,t} =  p_0^a + p^\delta_{0,t}										\label{eq:p:substation_active}	\\
				& p_{0,t} \geq -\mathbf{1}^\top\bp_t + \bp_t^\top{\mathbf{R}}\bp_t 
					+ \bq_t^\top{\mathbf{R}}\bq_t 										\label{eq:p:substation_losses}	\\
				& \mathbf{p}_t^{\top}\mathbf{f}_n\mathbf{f}_n^{\top}\mathbf{p}_t +  
					\mathbf{q}_t^{\top}\mathbf{f}_n\mathbf{f}_n^{\top} \mathbf{q}_t
					\leq \overline{S}_n,~\forall n \in \mathcal{N}						\label{eq:p:power_flows} 		\\
				& \underline{\bp}^d 	\leq \bp^{d} 	\leq \overline{\bp}^d			\label{eq:p:diesel_active} 		\\
				&  (p_{n,t}^r,\,q_{n,t}^r)\in \Omega_{n,t},~\forall n \in \mathcal{N}	\label{eq:p:inverter}			\\
				& \underline{v}_0 \leq v_0^a \leq \overline{v}_0 						\label{eq:p:voltage_substation} \\
				& \bv_t  = 2\mathbf{R}\bp_t + 2{\mathbf{X}}\bq_t + v_0^a \mathbf{1}		\label{eq:p:voltage_balance}	\\
				&  	\bv_t \in \mathcal{V}_B	 							\label{eq:p:voltage_hard}  		\\
				&	\E_t\left[\bv_t\right]\in \mathcal{V}_A				\label{eq:p:voltage_average}	\\
\mathrm{over} ~~& v_0^a, p_0^a, \bp^{d} ,  \{\bp_t, \bq_t, \bv_t, \bp^r_t, \bq^r_t, p_{0,t}, p^\delta_{0,t}\}_{t=1}^T.	\nonumber
\end{align}
\end{subequations}
The slow-timescale variables $\{v_0^a,p_0^a,\bp^{d}\}$ are set in advance, and remain fixed throughout the $T$ subsequent control slots over which the fast-timescale variables $\{\bp_t, \bq_t, \bv_t, \bp^r_t, \bq^r_t, p_{0,t}, p^\delta_{0,t}\}_{t=1}^T$ are implemented. The latter variables depend on the randomness of slot $t$ as well as slow-timescale decisions. 

Alternatively to \eqref{eq:p}, optimal grid operation can be posed as a \emph{probabilistic dispatch} that is identical to \eqref{eq:p} with the exception that \eqref{eq:p:voltage_average} is replaced by the probabilistic constraint
\begin{equation}\label{eq:prob}
\Pr\{\bv_t \notin \mathcal{V}_A\} \leq \alpha
\end{equation}
for some small parameter $\alpha>0$, say $\alpha=0.05$. The optimal cost for the probabilistic dispatch will be denoted by $\mathsf{P}_p^\ast$.

The objective function in \eqref{eq:p:objective} involves the cost of energy dispatched at the slow timescale plus the average fast-timescale energy management cost. Nodal (re)active power balance is ensured via \eqref{eq:p:active_balance}--\eqref{eq:p:reactive_balance}. Constraint \eqref{eq:p:substation_losses} accounts for the active power losses upon relaxing the quadratic equation in \eqref{eq:activelosses} to a convex inequality without loss of optimality. Constraint \eqref{eq:p:power_flows} limits the apparent power flow at each line $n$ based on \eqref{eq:appline}. Constraints \eqref{eq:p:voltage_substation}--\eqref{eq:p:voltage_average} are voltage regulation constraints: In detail, \eqref{eq:p:voltage_balance} relates squared voltage magnitudes to power injections [cf. \eqref{eq:voltage}]; \eqref{eq:p:voltage_substation} constraints the substation bus voltage; and \eqref{eq:p:voltage_hard} constraints voltages in $\mathcal{V}_B$. While \eqref{eq:p:voltage_average} maintains the average voltage magnitudes in $\mathcal{V}_A$, its alternative in \eqref{eq:prob} limits the probability of voltage magnitudes being outside $\mathcal{V}_A$.

\section{Problem Analysis} \label{s:analysis}
To facilitate algorithmic developments, the problem in \eqref{eq:p} is expressed in a compact form next. Collect the slow-timescale variables in vector $\bz^\top := [v_0^a, p_0^a, \bp^{d}]$; the fast-timescale variables at slot $t$ in $\by_t^\top := [\bp_t, \bq_t, \bv_t, \bp^r_t, \bq^r_t,p_{0,t},p_{0,t}^\delta]$; and the random variables involved at slot $t$ in $\bxi_t^{\top}:= [\overline{\mathbf{p}}_{t}^r, \bp_t^l, \bq_t^l]$. 

The constraints in \eqref{eq:p} can be classified into four groups:
\emph{(i)} Constraints involving fast-timescale variables only, such as \eqref{eq:p:reactive_balance}, \eqref{eq:p:power_flows}, \eqref{eq:p:inverter}, and \eqref{eq:p:voltage_hard}, that will be abstracted as $\by_t \in \mathcal{Y}_t$.\\
\emph{(ii)} Constraints \eqref{eq:p:diesel_active} and \eqref{eq:p:voltage_substation} that involve slow-timescale variables only, and they will be denoted as  $\bz \in \mathcal{Z}$.\\
\emph{(iii)} The \textit{linear} constraints \eqref{eq:p:active_balance}, \eqref{eq:p:substation_active}, and \eqref{eq:p:voltage_balance}, coupling slow- and fast-timescale variables as well as random variables. These constraints are collectively expressed as $\mathbf{K} \bz + \mathbf{B} \by_t = \mathbf{H}\bxi_t$ for appropriate matrices $\mathbf{K}$, $\mathbf{B}$, and $\mathbf{H}$.\\
\emph{(iv)} The ergodic constraints \eqref{eq:p:voltage_average} and \eqref{eq:prob} depend on the voltage sequence $\{\mathbf{v}_t\}_{t=1}^T$, hence coupling decisions across time. A substantial difference between \eqref{eq:p:voltage_average} and \eqref{eq:prob} is that the latter is a non-convex constraint.

Based on this grouping, the two dispatch problems can be compactly rewritten as
\begin{subequations}\label{eq:p_compact}
\begin{align}
\mathsf{P}_{(a,p)}^\ast := \min_{\bz, \{\by_t\}_{t=1}^T} 	~& f(\bz ) + \E_t\left[ g_t(\by_t )\right]		\label{eq:p1c:objective}\\
\text{s.to:}~~& \bz \in \mathcal{Z}															\label{eq:p1c:z}\\	
&\by_t \in \mathcal{Y}_t 								&\forall t 	\label{eq:p1c:y}\\		
&\mathbf{K} \bz + \mathbf{B} \by_t = \mathbf{H}\bxi_t	  &\forall t			\label{eq:p1c:linear}\\
&\E_t\left[\mathbf{h}(\by_t)\right]\leq \mathbf{0}\label{eq:p_compact_erg}
\end{align}
\end{subequations}
where $f(\bz ):=C_\textrm{D}(\bp^{d}) + \beta p_0^a$ and $g_t(\by_t):=C^t(p^\delta_{0,t}) + C_{\textrm{PV}}(\bp^r_t)$. For the average dispatch, the optimal cost in \eqref{eq:p_compact} is $\mathsf{P}_a^\ast$ and the function in \eqref{eq:p_compact_erg} is $\mathbf{h}(\by_t)=[\mathbf{v}_t-\overline{v}_A\mathbf{1}, \underline{v}_A\mathbf{1}-\mathbf{v}_t]$. For the probabilistic dispatch, the optimal cost is $\mathsf{P}_p^\ast$ and the function in \eqref{eq:p_compact_erg} is $\mathbf{h}(\by_t)=\mathbbm{1}\{\mathbf{v}_t \notin \mathcal{V}_A\} - \alpha$. 

The optimal values for the slow-timescale variables $\bz$ 
must 
be decided in advance. Once the optimal $\bz$ 
is 
found, it remains fixed over the slow-timescale interval. The fast-timescale decisions $\by_t(\bz)$ for slot $t$ depend on 
$\bz$, while the subscript $t$ indicates their dependence on the realization $\bxi_t$. Both the average and the probabilistic dispatch are stochastic programming problems with recourse~\cite{VaWuBi11}. Their costs can be decomposed as $\mathsf{P}_{(a,p)}^\ast = \min_{\bz\in \mathcal{Z}} f(\bz ) + G_{(a,p)}(\bz)$, where the so termed \emph{expected recourse function} is defined as
\begin{subequations}\label{eq:G1_recourse}
\begin{align}
G_{(a,p)}(\bz) := \min_{\{\by_t\in \mathcal{Y}_t\}} 	~& \E_t\left[ g_t(\by_t )\right]		\label{eq:Gr:objective_avg}\\
\text{s.to:} ~~ 	&\mathbf{K} \bz + \mathbf{B} \by_t = \mathbf{H}\bxi_t		 & \forall t	\label{eq:Gr:linear_avg}\\
&\E_t\left[\mathbf{h}(\by_t)\right]\leq \mathbf{0}.\label{eq:G1r:erg}
\end{align}
\end{subequations}
Since problem \eqref{eq:G1_recourse} depends on $\bz$, its minimizer can be written as $\{\by_t^*(\bz)\}_{t=1}^T$ and the recourse function as $G_{(a,p)}(\bz)=\E_t[ g_t(\by_t^*(\bz) )]$. The ensuing two sections solve the average and the probabilistic dispatches.

\section{Average Dispatch Algorithm}\label{s:ada}
This section tackles problem \eqref{eq:p_compact} with the ergodic constraint in \eqref{eq:p_compact_erg}, for which $\mathbf{h}(\by_t)=[\underline{v}_A\mathbf{1}-\mathbf{v}_t, \mathbf{v}_t-\overline{v}_A\mathbf{1}]$. Although convex, problem \eqref{eq:p_compact} is challenging due to the coupling across $\{\mathbf{y}_t\}_{t=1}^T$ and between fast- and slow-timescale variables. Dual decomposition is adopted to resolve the coupling across $\{\by_t\}_{t=1}^T$. The partial Lagrangian function for \eqref{eq:G1_recourse} is
\begin{align}\label{eq:Lagrangian1}
L_a(\{\by_t\},\bnu) := \E_t\left[g_t(\by_t)+ \bnu^\top\mathbf{h}(\by_t) \right]
\end{align}
with the entries of $\bnu$ being the multipliers associated with the upper and lower per-bus constraints in \eqref{eq:p_compact_erg}. The corresponding dual function is
\begin{align}
D_a(\bnu;\bz) := \min_{\{\by_t\in \mathcal{Y}_t\}}~&L_a(\{\by_t\},\bnu)\label{eq:dualf}\\
\text{s.to:}~&\mathbf{K} \bz + \mathbf{B} \by_t = \mathbf{H}\bxi_t \quad \forall t.\nonumber
\end{align}
Observe that after dualizing, the minimization in \eqref{eq:dualf} is separable over the realizations $\{\bxi_t\}$. Precisely, the optimal fast-timescale variable for fixed $(\bnu,\bz)$ and for a specific realization $\bxi_t$ can be found by solving:
\begin{subequations}\label{eqs:inst_opt1}
\begin{align}
\by_t^*(\bnu,\bz)\in\arg\min_{\by_t\in \mathcal{Y}_t} ~&g_t(\by_t) + \bnu^\top\mathbf{h}(\by_t)\label{eqs:inst_opt1:cost}\\
	\text{s.to:} ~&\mathbf{K} \bz + \mathbf{B} \by_t = \mathbf{H}\bxi_t.\label{eqs:inst_opt1:con}
\end{align}
\end{subequations}
For future reference, let us also define $\blambda_t^*(\bnu,\bz)$ as the optimal Lagrange multiplier associated with \eqref{eqs:inst_opt1:con}. If $\bnu$ is partitioned as $\bnu^\top=[\underline{\bnu}^\top, \overline{\bnu}^\top]$ with $\underline{\bnu}$ corresponding to constraint $\mathbb{E}_t[\bv_t]\geq \underline{v}_A\mathbf{1}$ and $\overline{\bnu}$ to $\mathbb{E}_t[\bv_t]\leq \overline{v}_A\mathbf{1}$, then \eqref{eqs:inst_opt1} simplifies to
\begin{align}\label{eq:expansion_instantaneous}
\by_t^*(\bnu,\bz)\in\argmin &~ C^t(p^\delta_{0,t}) + C_{\textrm{PV}}(\bp^r_t)+ (\overline{\bnu}-\underline{\bnu})^{\top}\bv_t\\
\text{s.to:} &~ \eqref{eq:p:active_balance}-\eqref{eq:p:power_flows}, 
\eqref{eq:p:inverter}, \eqref{eq:p:voltage_balance}, \eqref{eq:p:voltage_hard}  \nonumber\\
\mathrm{over}&~ \{\bp_t, \bq_t, \bv_t, \bp^r_t, \bq^r_t, p_{0,t},p_{0,t}^\delta\}	\nonumber
\end{align}
and can be solved as a convex quadratic program. Given the optimal pair $(\bnu^*,\bz^*)$, the optimal fast-timescale variables $\by_t$ can be thus found for any $\bxi_t$. 

Back to finding the optimal primal and dual slow-timescale variables, note that the dual problem associated with \eqref{eq:dualf} is
\begin{equation}\label{eq:dualp}
\bnu^*:=\arg\max_{\bnu\geq \mathbf{0}}~D_a(\bnu; \bz).
\end{equation}
Duality theory asserts that \eqref{eq:dualp} is a convex problem. Moreover, assuming a strictly feasible point exists for \eqref{eq:G1_recourse}, strong duality implies that $G_a(\bz) = D_a(\bnu^*,\bz)$. Due to the latter, the original problem in \eqref{eq:p_compact} can be transformed to:
\begin{subequations}\label{eq:ftilde}
\begin{align}
\min_{\bz \in \mathcal{Z}} f(\bz) + G_a(\bz) &= \min_{\bz \in \mathcal{Z}} \{f(\bz) + \max_{\bnu\geq \mathbf{0}} D_a(\bnu;\bz)\}\label{eq:ftilde1}\\
&=\min_{\bz \in \mathcal{Z}} \max_{\bnu \geq \mathbf{0}} \tilde{f}_a(\bnu, \bz)\label{eq:ftilde2}
\end{align}
\end{subequations}
where the auxiliary function $\tilde{f}_a$ is defined as:
\begin{equation}
\tilde{f}_a(\bnu, \bz) := f(\bz) + D_a(\bnu;\bz).
\end{equation}
Being a dual function, $D_a(\bnu;\bz)$ is a concave function of $\bnu$. At the same time, $D_a(\bnu;\bz)$ is a perturbation function with respect to $\bz$; and hence, it is a convex function of $\bz$~\cite{BoVa04}. Recall that $f(\bz)$ is a convex function of $\bz$ too. Therefore, the auxiliary function $\tilde{f}_a(\bnu, \bz)$ is convex in $\bz$ and concave in $\bnu$. Because of the randomness of $\{\bxi_t\}$, function $D_a(\bnu;\bz)$ in \eqref{eq:dualf} is stochastic. Consequently, problem \eqref{eq:ftilde2} is a stochastic convex-concave saddle point problem \cite{BoVa04}, \cite{Nemirovski09}. 

To solve \eqref{eq:ftilde2}, we rely on the stochastic saddle-point approximation method of~\cite{Nemirovski09}. The method involves the subgradient of $\tilde{f}_a$ with respect to $\bz$, and its supergradient with respect to $\bnu$. Upon viewing $D_a(\bnu, \bz)$ in \eqref{eq:dualf} as a perturbation function of $\bz$, the subgradient of $\tilde{f}_a$ with respect to $\bz$ is~\cite{BoVa04}
\begin{equation}\label{eq:g_z}
\partial_{\bz} \tilde{f}_a = \partial_{\bz} f(\bz) + \mathbf{K}^\top\E_t[\blambda_t^*(\bnu,\bz)].
\end{equation}
By definition of the dual function, the supergradient of $\tilde{f}_a$ with respect to $\bnu$ is
\begin{equation}\label{eq:g_nu}
\partial_{\bnu}\tilde{f}_a = \E_t[\mathbf{h}(\by_t^*(\bnu,\bz))].
\end{equation}
The stochastic saddle point approximation method of \cite{Nemirovski09} involves primal-dual subgradient iterates with the expectations in \eqref{eq:g_z}--\eqref{eq:g_nu} being replaced by their instantaneous estimates based on a single realization $\bxi_k$. Precisely, the method involves the iterates over $k$:
\begin{subequations}\label{eqs:primal_dual1}
\begin{align}
\bnu^{k+1} &:= [\bnu^k + \diag(\bmu_k)\mathbf{h}(\by_k^*(\bnu^k,\bz^k))]_+\\
\bz^{k+1} &:= [\bz^k - \diag(\boldsymbol{\epsilon}_k) (\partial_z f(\bz^k) + \mathbf{K}^\top\boldsymbol{\lambda}_k^*(\bnu^k,\bz^k))]_\mathcal{Z}\label{eq:primal_dual_z}
\end{align}
\end{subequations}
where the operator $[\cdot]_{\mathcal{Z}}$ projects its argument onto $\mathcal{Z}$; and vectors $\bmu_k=\bmu_0/\sqrt{k}$ and $\bepsilon_k=\bepsilon_0/\sqrt{k}$ collect respectively the primal and dual step sizes for positive $\bmu_0$ and $\bepsilon_0$. At every iteration $k$, the method draws a realization $\bxi_k$ and solves \eqref{eqs:inst_opt1} for the tuple $(\bxi_k,\bnu^k,\bz^k)$ to acquire $(\by_k^*(\bnu^k,\bz^k),\blambda_k^*(\bnu^k,\bz^k))$ and perform the primal-dual updates in \eqref{eqs:primal_dual1}. The method finally outputs the sliding averages of the updates as:
\begin{subequations}\label{eq:slide}
\begin{align}
\tilde{\bz}^k &:={\textstyle \big(\sum_{i=\lceil k/2\rceil}^k \bz^i /\sqrt{i}\big)/\big(\sum_{i=\lceil k/2\rceil}^k 1/\sqrt{i}\big)}\\
\tilde{\bnu}^k &:={\textstyle \big(\sum_{i=\lceil k/2\rceil}^k \bnu^i /\sqrt{i}\big)/\big(\sum_{i=\lceil k/2\rceil}^k 1/\sqrt{i}\big).}
\end{align}
\end{subequations}
The proposed scheme converges to the value $\tilde{f}_a(\bnu^*, \bz^*)$ obtained at a saddle point $(\bnu^*, \bz^*)$ asymptotically in the number of iterations $k$~\cite[Sec.~3.1]{Nemirovski09}.

\begin{algorithm}[t]
	\caption{Average Dispatch Algorithm (ADA)}
	\label{alg:avg}
	\begin{algorithmic}[1]
		\State Initialize $(\bz^0,\bnu^0$).
		\Repeat ~for $k = 0,1,\ldots$
		\State Draw sample $\bxi_k$.
		\State Find $(\by_k^*(\bnu^k,\bz^k),\blambda_k^*(\bnu^k,\bz^k))$ by solving \eqref{eqs:inst_opt1}.
		\State Update $(\bz^{k+1},\bnu^{k+1})$ via \eqref{eqs:primal_dual1}.
		\State Compute sliding averages $(\tilde{\bz}^k, \tilde{\bnu}^k)$ through \eqref{eq:slide}.
		\Until convergence of $(\tilde{\bz}^k, \tilde{\bnu}^k)$.
		\State \textbf{Output} $ \bz^{*} = \tilde{\bz}^k$ and $\bnu^{*} = \tilde{\bnu}^k$.
	\end{algorithmic}
\end{algorithm}

Upon convergence of the iterates in \eqref{eq:slide}, the slow-timescale variables $\bz^*$ have been derived together with the optimal Lagrange multiplier $\bnu^*$ related to constraint \eqref{eq:G1r:erg}. The grid operator can implement  $\bz^*$, and the fast-timescale decisions $\by_t^*$ for a realization $\bxi_t$ can be found by solving \eqref{eq:expansion_instantaneous}. The average dispatch algorithm (ADA) is summarized as Alg.~\ref{alg:avg}.

\section{Probabilistic Dispatch Algorithm}\label{s:pda}
The probabilistic version of problem \eqref{eq:p_compact} is considered next. Here, the ergodic constraint \eqref{eq:p_compact_erg} reads $h(\by_t)= \mathbbm{1}\{\mathbf{v}_t \notin \mathcal{V}_A\} - \alpha$. Despite the non-convexity of the probabilistic constraint, \eqref{eq:G1_recourse} can still be solved optimally. However, optimality for \eqref{eq:p_compact} cannot be guaranteed. A heuristic solution is detailed next by adapting the solution of Sec.~\ref{s:ada}.

To that end, dual decomposition is used here as well. If $\nu$ is the scalar Lagrange multiplier associated with constraint \eqref{eq:p_compact_erg}, the partial Lagrangian function for \eqref{eq:G1_recourse} is now $L_p(\{\by_t\},\nu):= \E_t\left[g_t(\by_t)+ \nu(\mathbbm{1}\{\bv_t \notin \mathcal{V}_A\} - \alpha) \right]$. The corresponding dual function, fast-timescale problem, and dual problem are defined analogously to \eqref{eq:dualf}, \eqref{eqs:inst_opt1}, and \eqref{eq:dualp}. The indicator function renders $L_p(\{\by_t\},\nu)$ non-convex. Surprisingly enough though, under the practical assumption that $\{\bxi_t\}$ follows a continuous pdf, problem \eqref{eq:G1_recourse} enjoys zero duality gap; see~\cite[Th.~1]{ribeiro10}.

The additional challenge here is the non-convexity of the Lagrangian minimization:
\begin{align}\label{eqs:inst_opt2}
\by_t^*(\nu,\bz)\in\arg\min_{\by_t\in \mathcal{Y}_t} ~&g_t(\by_t)+ \nu \mathbbm{1}\{\bv_t \notin \mathcal{V}_A\}\\
\text{s.to:}~~&\mathbf{K} \bz + \mathbf{B} \by_t = \mathbf{H}\bxi_t.\nonumber
\end{align}
Because the indicator function takes only the values $\{0,1\}$ however, the solution to \eqref{eqs:inst_opt2} can be found by solving a pair of slightly different convex problems. The first problem is 
\begin{subequations}\label{eq:J_tight}
		\begin{align}
		\by_{t,A}^*(\bz)\in\arg\min_{\by_t \in  \mathcal{Y}_t} ~&~ g_t(\by_t)\\
		\text{s.to:} ~~&~\mathbf{K} \bz + \mathbf{B} \by_t = \mathbf{H}\bxi_t\label{eq:J_tight_con}\\
		&~\bv_t\in\mathcal{V}_A
		\end{align}
	\end{subequations}
whereas the second problem ignores constraint $\mathbf{v}_t \in \mathcal{V}_A$ as
	\begin{subequations}\label{eq:J_loose}
		\begin{align}
		\by_{t,B}^*(\bz)\in \arg \min_{\by_t \in \mathcal{Y}_t} ~&~ g_t(\by_t)\\
		\text{s.to:} ~~&~\mathbf{K} \bz + \mathbf{B} \by_t = \mathbf{H}\bxi_t.\label{eq:J_loose_con}
		\end{align}
	\end{subequations}
From the point of view of \eqref{eqs:inst_opt2}, if the voltages in $\by_{t,B}^*(\bz)$ do not belong to $\mathcal{V}_A$, the solution to the second problem will incur an additional cost quantified by $\nu$. Observe that neither problem \eqref{eq:J_tight} nor \eqref{eq:J_loose} depend on $\nu$, while their complexity is similar to the one problem \eqref{eqs:inst_opt1}. Suppose that \eqref{eq:J_tight} and \eqref{eq:J_loose} have been solved and let $\blambda_{t,A}^*(\bz)$ and $\blambda_{t,B}^*(\bz)$ denote the optimal multipliers associated with  \eqref{eq:J_tight_con} and \eqref{eq:J_loose_con}, respectively. Then, problem \eqref{eqs:inst_opt2} can be neatly tackled by identifying two cases:

\begin{algorithm}[t]
	\caption{Probabilistic Dispatch Algorithm (PDA)}
	\label{alg:prob2}
	\begin{algorithmic}[1]
		\State Initialize $(\bz^0,\nu^0)$.
		\Repeat~for $k = 0,1,\ldots$
		\State Draw sample $\bxi_k$.
		\State Find $(\by_{k,B}^*(\nu^k,\bz^k),\blambda_{k, B}^*(\nu^k,\bz^k))$ by solving \eqref{eq:J_loose}.
		\State Set $\by_t^*(\nu, \bz):=\by_{t,B}^*(\bz)$ and $\blambda_t^*(\nu, \bz):=\blambda_{t,B}^*(\bz)$.
		\If{$\bv_{k,B}^*(\bz) \notin \mathcal{V}_A$,}
		{find $\by_{k,A}^*(\nu^k,\bz^k)$ and $\blambda_{k, A}^*(\nu^k,\bz^k)$ by solving \eqref{eq:J_tight}.}
		{\If{$g_t(\by_{t,A}^*(\bz)) \leq g_t(\by_{t,B}^*(\bz)) +\nu$,}
			{set $\by_t^*(\nu, \bz):=\by_{t,A}^*(\bz)$ and $\blambda_t^*(\nu, \bz):=\blambda_{t,A}^*(\bz)$.}
			\EndIf}
		\EndIf
		\State Update $(\bz^{k+1},\nu^{k+1})$ via \eqref{eqs:primal_dual2}.
		\State Compute sliding averages $(\tilde{\bz}^k, \tilde{\nu}^k)$ through \eqref{eq:slide}.
		\Until convergence of $(\tilde{\bz}^k, \tilde{\nu}^k)$.
		\State \textbf{Output} $ \bz^{*} = \tilde{\bz}^k$ and $\nu^{*} = \tilde{\nu}^k$.
	\end{algorithmic}
\end{algorithm}

\textbf{(c1)} If $g_t(\by_{t,A}^*(\bz))> g_t(\by_{t,B}^*(\bz)) +\nu$, then $\by_{t,B}^*(\bz)$ is a minimizer of \eqref{eqs:inst_opt2} as well and voltages are allowed to lie outside $\mathcal{V}_A$. In this case, set $\by_t^*(\nu, \bz):=\by_{t,B}^*(\bz)$ and $\blambda_t^*(\nu, \bz):=\blambda_{t,B}^*(\bz)$. This case includes instances where problem \eqref{eq:J_tight} is infeasible for which $g_t(\by_{t,A}^*(\bz))=\infty$.

\textbf{(c2)} If $g_t(\by_{t,A}^*(\bz))\leq g_t(\by_{t,B}^*(\bz)) +\nu$, then $\by_{t,A}^*(\bz)$ minimizes \eqref{eqs:inst_opt2} too and voltages lie within $\mathcal{V}_A$. In this case, set $\by_t^*(\nu, \bz) := \by_{t,A}^*(\bz)$ and $\blambda_t^*(\nu, \bz):=\blambda_{t,A}^*(\bz)$. 

Case (c2) covers also instances where $\bv_{t,B}^*(\bz)$ happens to lie in $\mathcal{V}_A$. In these particular instances, $\by_{t,B}^*(\bz)$ serves as a minimizer of \eqref{eq:J_tight} too. Then, it follows that $g_t(\by_{t,A}^*(\bz))=g_t(\by_{t,B}^*(\bz))\leq g_t(\by_{t,B}^*(\bz)) +\nu$ for $\nu\geq 0$. This implies that one can solve \eqref{eq:J_loose} first and, if $\bv_{t,B}^*(\bz)\in\mathcal{V}_A$, there is no need to solve problem \eqref{eq:J_tight}.

To find the optimal slow-timescale variables under the probabilistic dispatch, the stochastic primal-dual iterations of Sec.~\ref{s:ada} are adapted here as 
\begin{subequations}\label{eqs:primal_dual2}
	\begin{align}
	\nu^{k+1} &:= [\nu^k + \mu_k(\mathbbm{1}\{\bv^*_k(\nu^k, \bz^k) \notin \mathcal{V}_A\} -\alpha)]_+ \\
	\bz^{k+1} &:= [\bz^k - \diag(\boldsymbol{\epsilon}_k) (\partial_z f(\bz^k) + \mathbf{K}^\top \blambda_k^*(\nu^k, \bz^k)]_\mathcal{Z}.\label{eqs:primal_dual2:z}
	\end{align}
\end{subequations}
The probabilistic dispatch algorithm (PDA) is tabulated as Alg.~\ref{alg:prob2}. Because function $G_p(\bz)$ is not necessarily convex, the iterates in \eqref{eqs:primal_dual2} are not guaranteed to converge to a minimizer of \eqref{eq:p_compact}. The practical performance of PDA in finding $\bz^\ast$ is numerically validated in Sec.~\ref{s:numerical}.

\section{Numerical Tests}\label{s:numerical}
The proposed grid dispatches were tested on a 56-bus Southern California Edison (SCE) distribution feeder~\cite{tac2015gan}. 5-MW PVs were added on buses 44 and 50; both with 6-MVA inverters enabling power factors as low as 0.83 (leading or lagging) at full solar generation. The prices for the energy exchange with the main grid were $\beta = 37$~\$/MWh; $\gamma_b = 45$~\$/MWh, and $\gamma_s = 19$~\$/MWh. Diesel generators with capacity $\overline{p}_n^d = 0.5$ MW were sited on buses 10, 18, 21, 30, 36, 43, 51, and 55. The cost of diesel generation was $C_\textrm{D} (\bp^d)= \sum_{n =1}^N(30 p_n^d + 15 (p_n^d)^2)$~\$/h with $\bp^d$ expressed in MW. Apparent power flows were limited to 7 MVA. The voltage operation limits were set to $\underline{v}_{A}=0.98^2$, $\overline{v}_{A}=1.02^2$, $\underline{v}_{B}=0.97^2$, and $\overline{v}_{B}=1.03^2$, expressed in pu with respect to a voltage base of 12~kV. (Re)active nodal loads were Gaussian distributed with the nominal load of the SCE benchmark as mean value, and standard deviation of 0.2 times the nominal load. The solar energy generated at each PV was drawn uniformly between 0.5 and 1 times the actual power PV rating.

\begin{figure}[t]
\centering
\includegraphics[width=0.97\columnwidth]{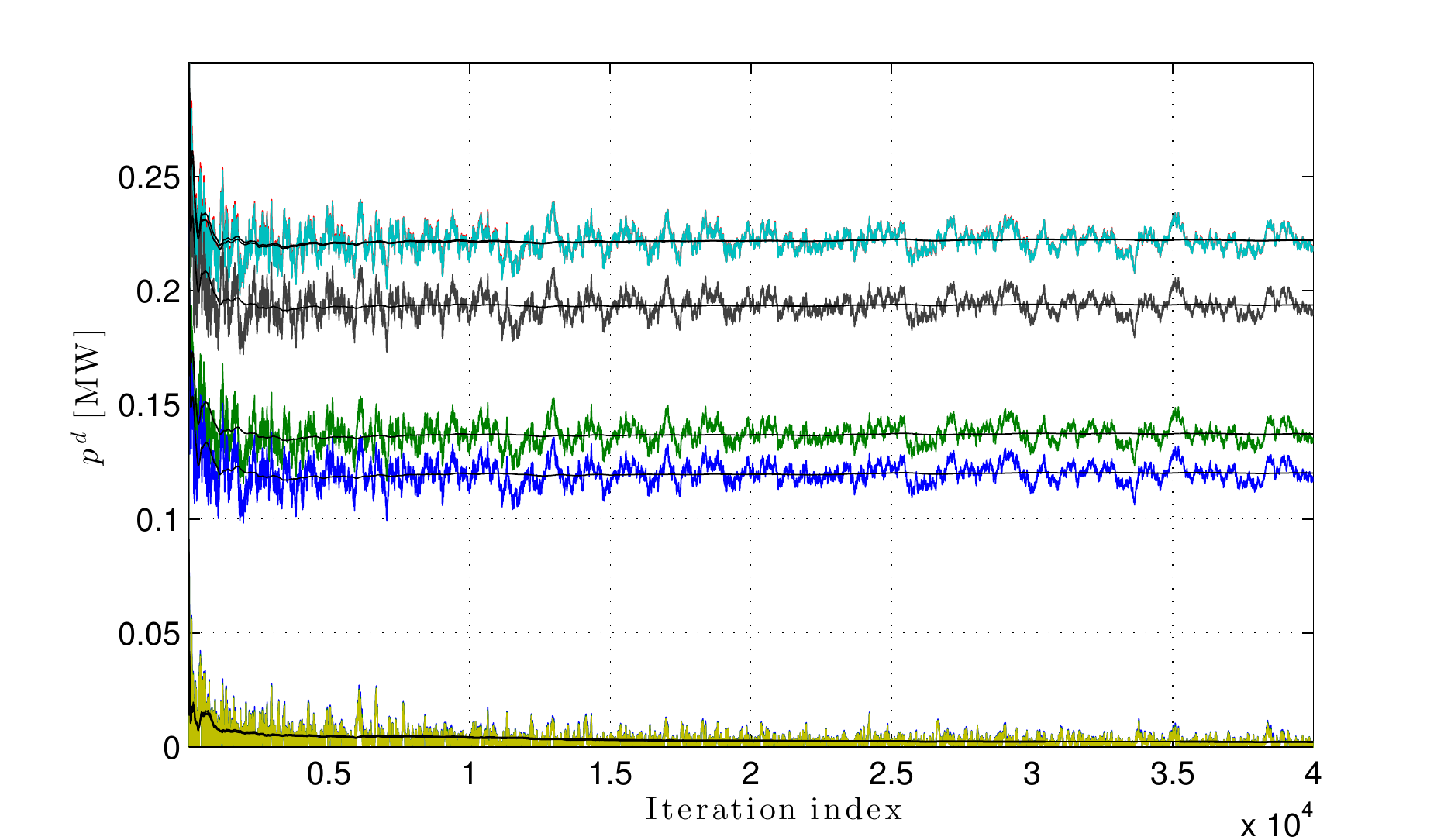}
\includegraphics[width=0.97\columnwidth]{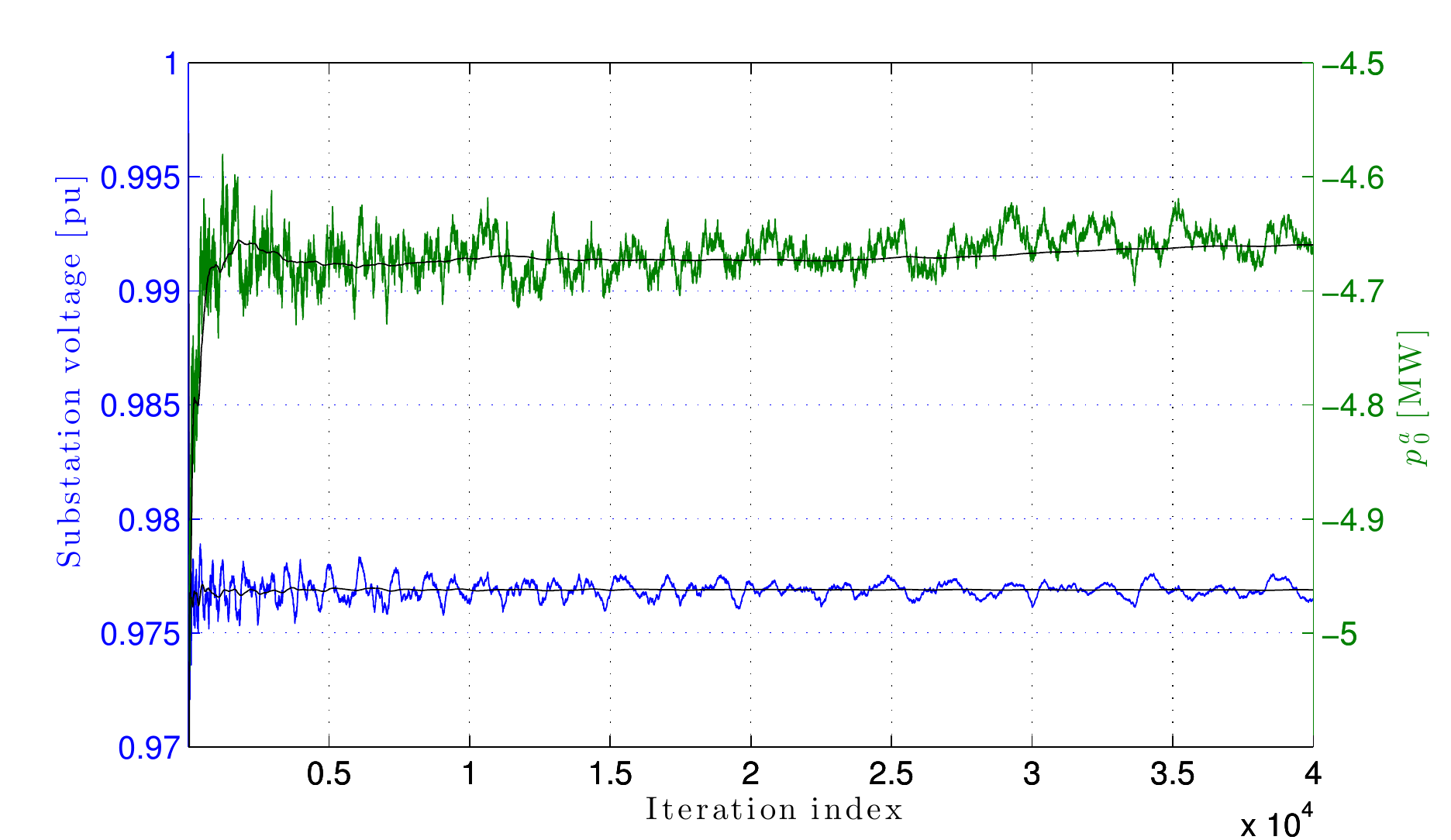}
\caption{Convergence of primal variables for ADA: (top) diesel generation; (bottom) substation voltage $v_0$ (left y-axis) and energy exchange $p_0^a$ (right y-axis). Sliding averages of optimization variables are depicted too.}
\vspace*{-1em}
\label{f:convergence_avg}
\end{figure}

\begin{figure}[t]
\centering
\includegraphics[width=\columnwidth]{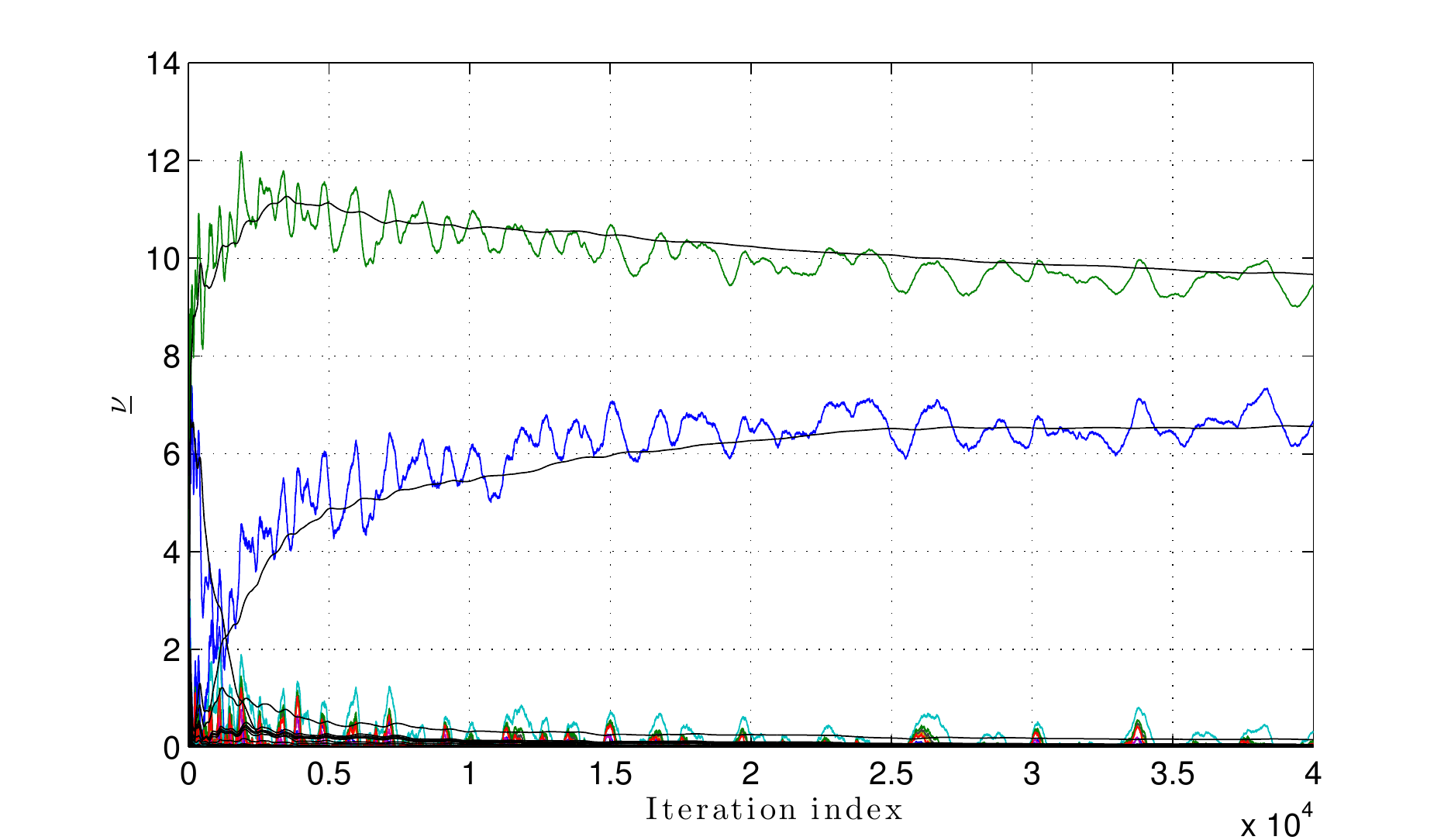}
\includegraphics[width=\columnwidth]{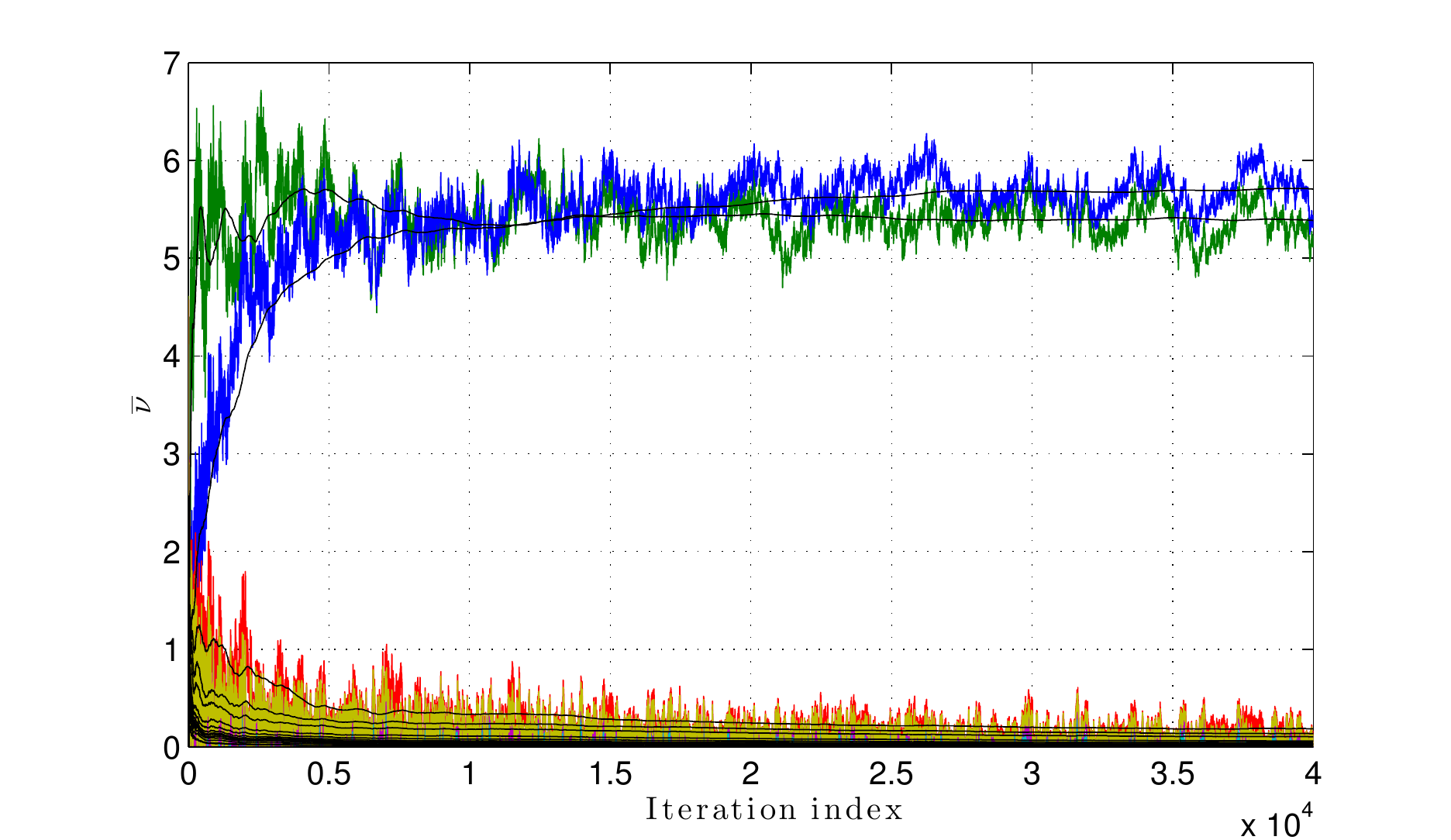}
\caption{Convergence of dual variables for ADA: (top) dual variables associated with average lower voltage limits for all buses; and (bottom) dual variables associated with average upper voltage limits for all buses. Sliding averages of optimization variables are depicted too.}\vspace{-5mm}
\label{f:convergence_avg}
\end{figure}

ADA was run 
with step sizes proportional to $1/\sqrt{k}$ with initial values $\epsilon_0^{v_0}=4\cdot 10^{-5}$, $\epsilon_0^{p_0}=4\cdot 10^{-1}$, $\epsilon_0^{p_d}=6\cdot 10^{-3}$, and $\mu_0=225$, to account for different dynamic ranges. The iterates for primal and dual variables as well as their corresponding sliding averages are depicted in Fig.~\ref{f:convergence_avg}. Primal and dual slow-timescale variables hover in a small range whose width diminishes with time. Their sliding averages converge asymptotically. The algorithm reaches a practically meaningful solution within 5,000 iterations. Buses 44 and 50 are prone to overvoltages since they host PV generation, 
and buses 2 and 15 are prone to under-voltages; 
thus yielding non-zero dual variables for the average upper and lower voltage constraints, respectively.

\begin{figure*}[t] 
\centering
\includegraphics[width=0.32\textwidth]{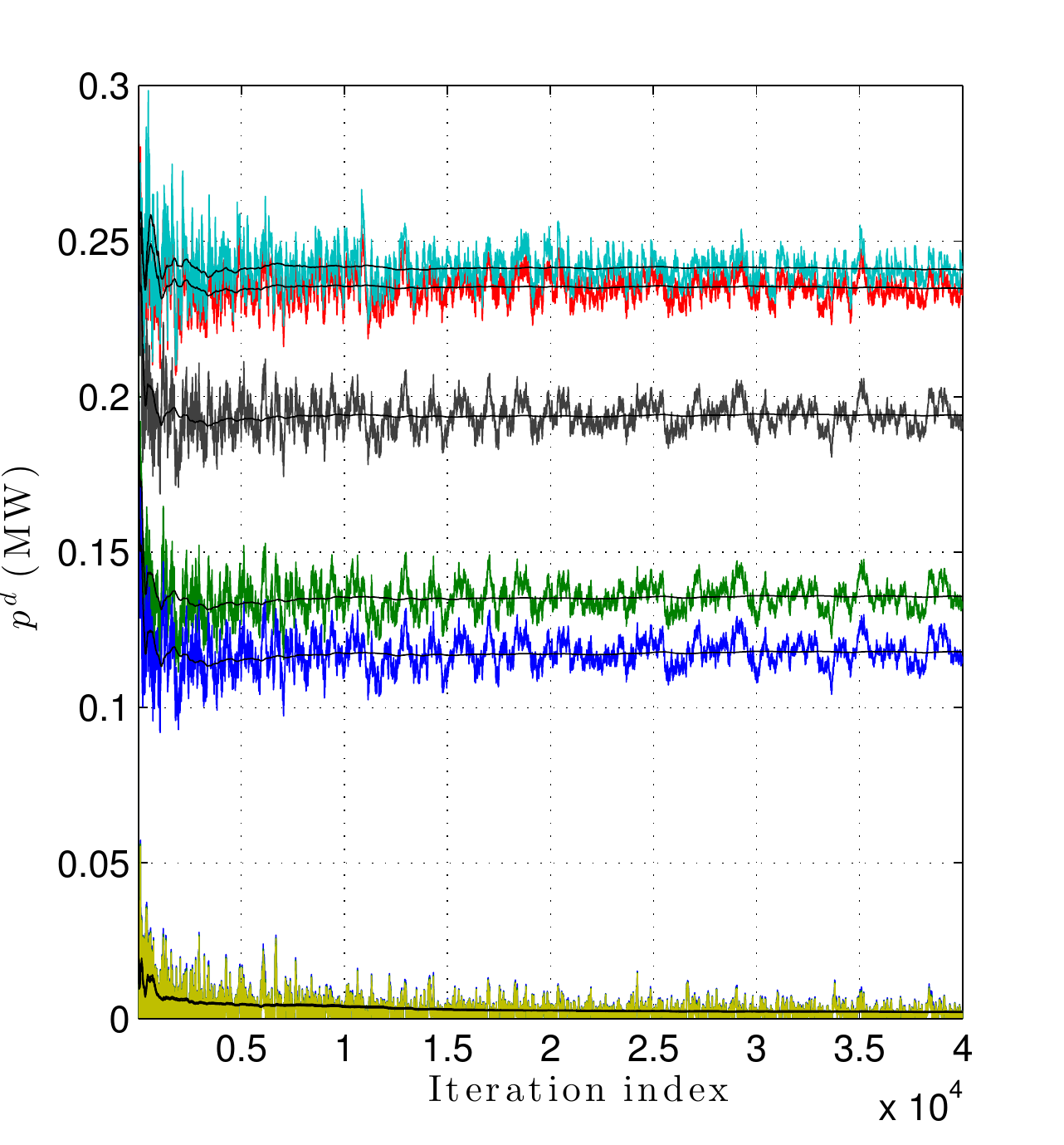}
\includegraphics[width=0.32\textwidth]{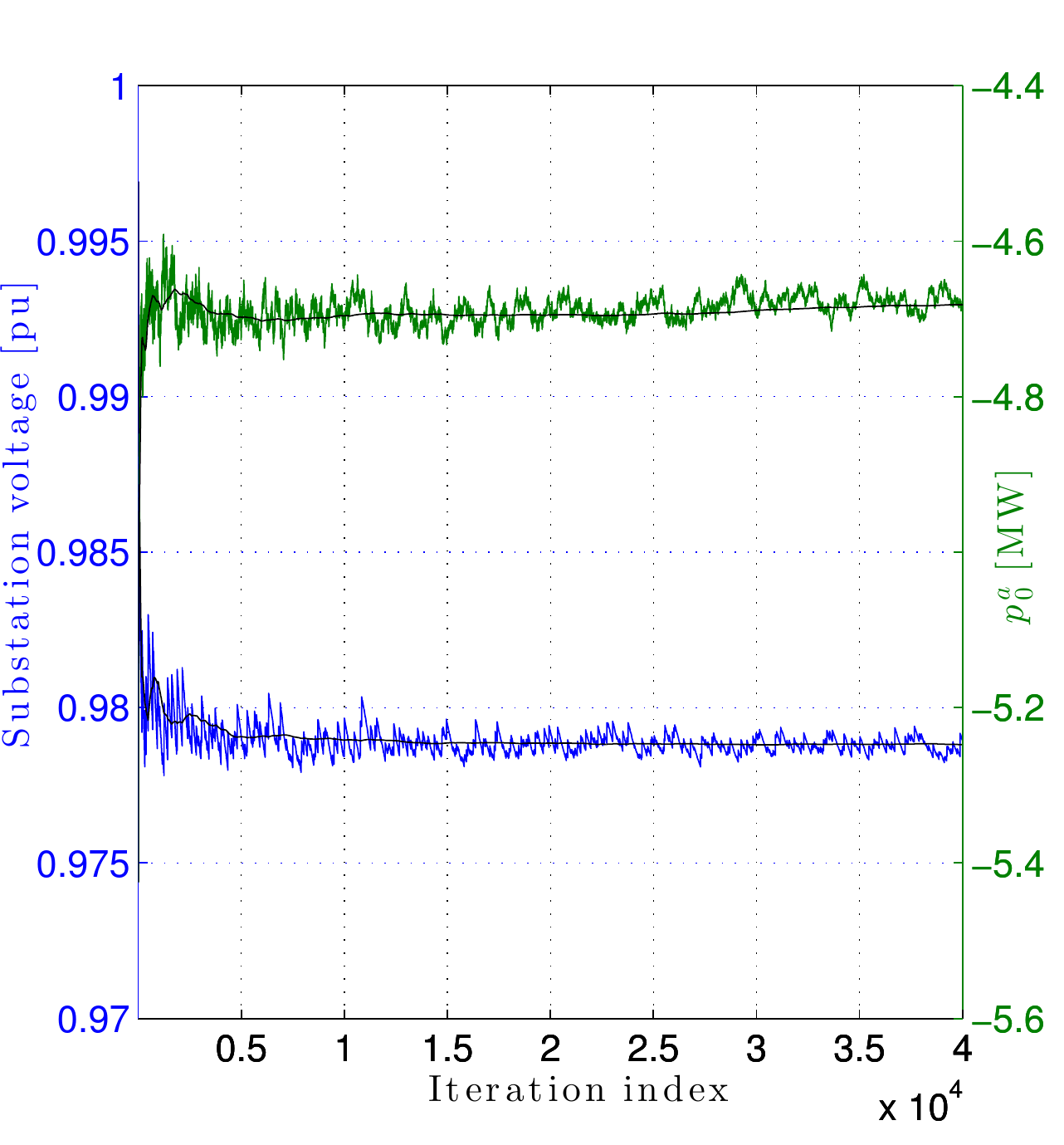}
\includegraphics[width=0.32\textwidth]{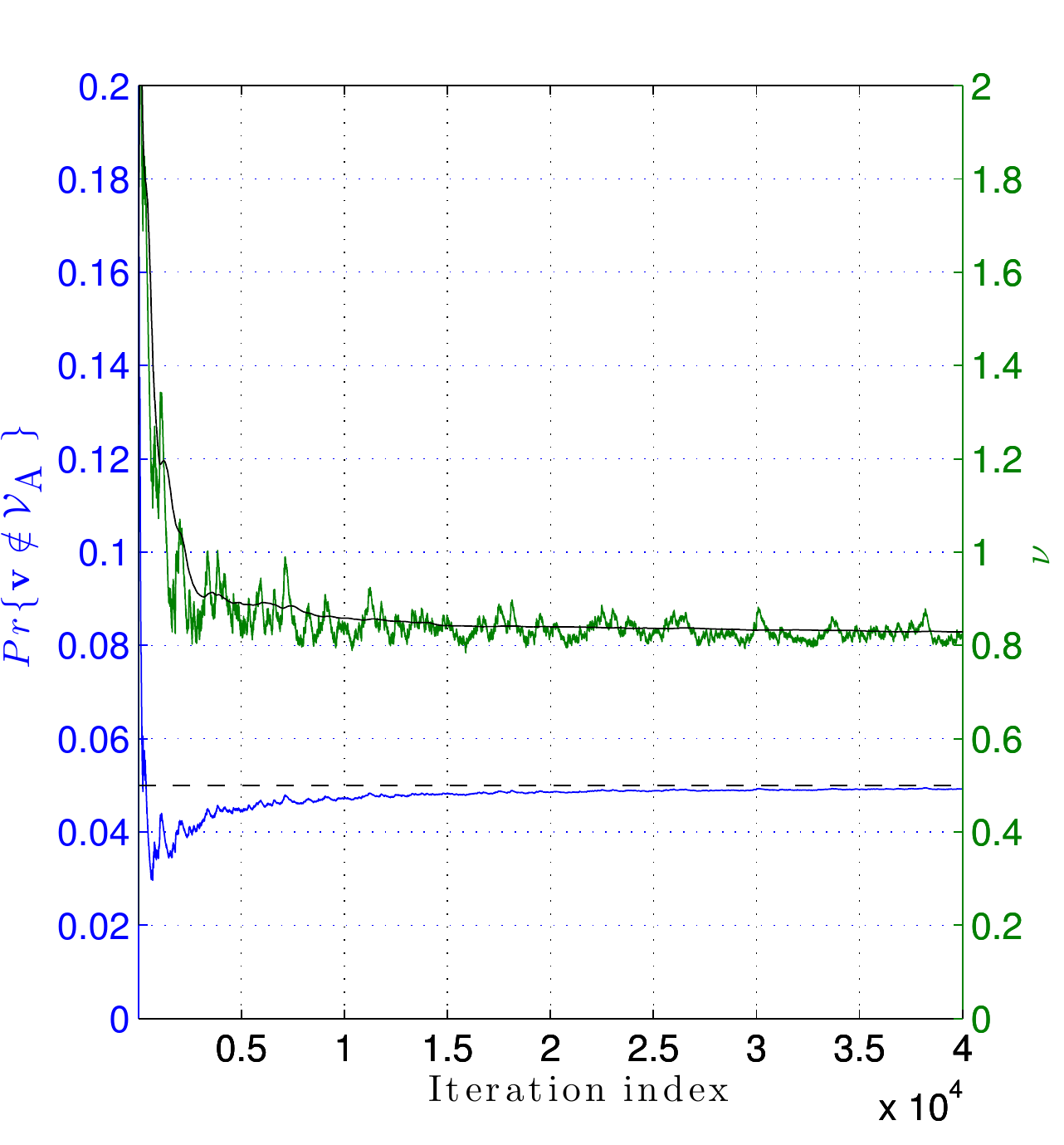}
\vspace*{-1em}
\caption{Convergence for PDA: (left) diesel generation; (middle) substation voltage (left y-axis) and energy exchange $p_0^a$ (right y-axis); and (right) dual variable related to probabilistic constraint (left y-axis) and under-/over-voltage probability (right y-axis). Sliding averages of optimization variables are shown too.}
\label{f:convergence_prob}
\end{figure*}

\begin{figure}[t]
\centering
\includegraphics[width=\columnwidth]{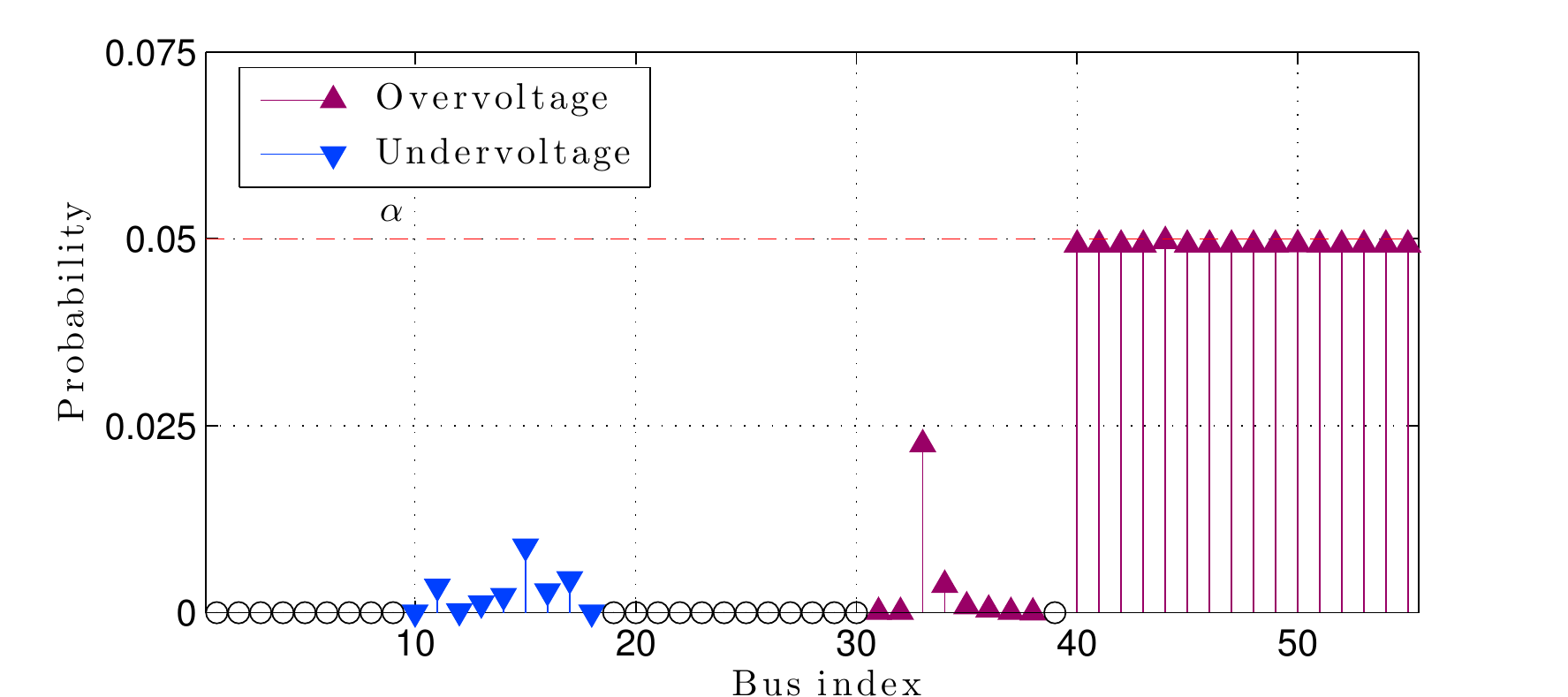}
\caption{Per-bus probability of under-/over-voltages.}
\label{f:probabilities}
\end{figure}

PDA was tested using the same simulation setup for $\alpha= 0.05$ and $\mu_0=1$. Figure~\ref{f:convergence_prob} shows the convergence of primal and dual variables, and the probability of voltages deviating from $\mathcal{V}_A$. Granted that the probabilistic constraint in \eqref{eq:prob} applies collectively to all buses, the under-/over-voltage probabilities on a per-bus basis is depicted in Fig.~\ref{f:probabilities}. The occurrences of overvoltage seem to be shared primarily among buses 40--56 which are neighboring to the PV buses 40 and 55. On the contrary, buses 10--16 being electrically far from both the substation and PVs, experience under-voltage with a small probability.

\begin{figure}[t]
\centering
\includegraphics[width=\columnwidth]{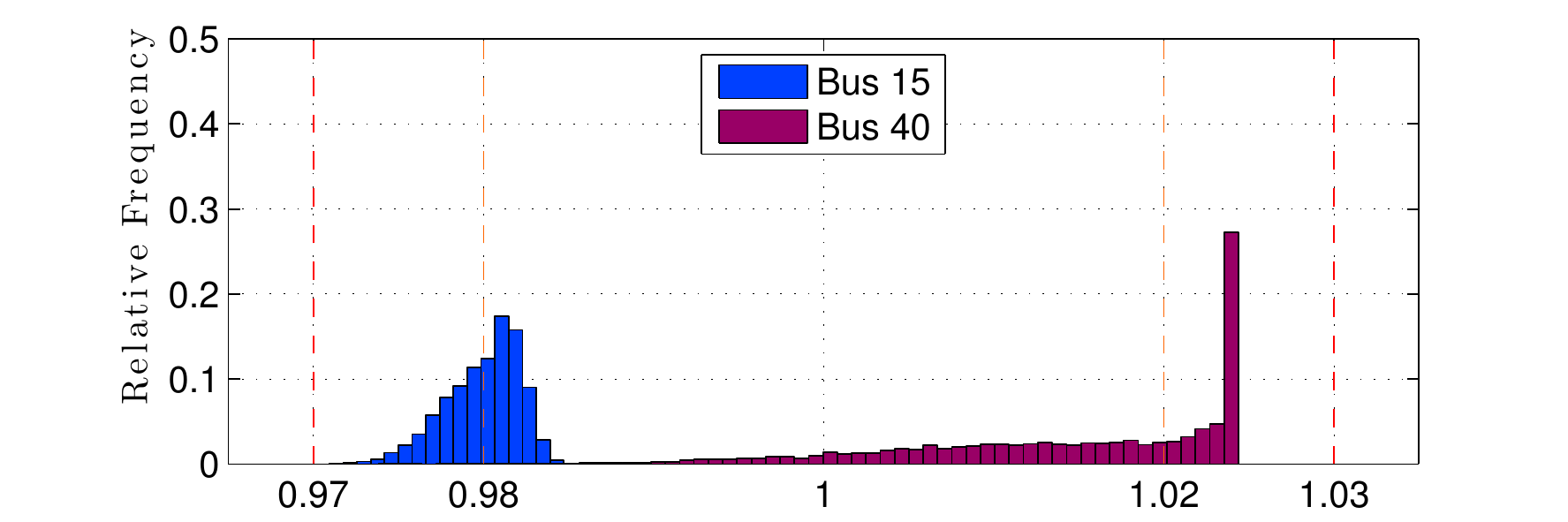}
\includegraphics[width=\columnwidth]{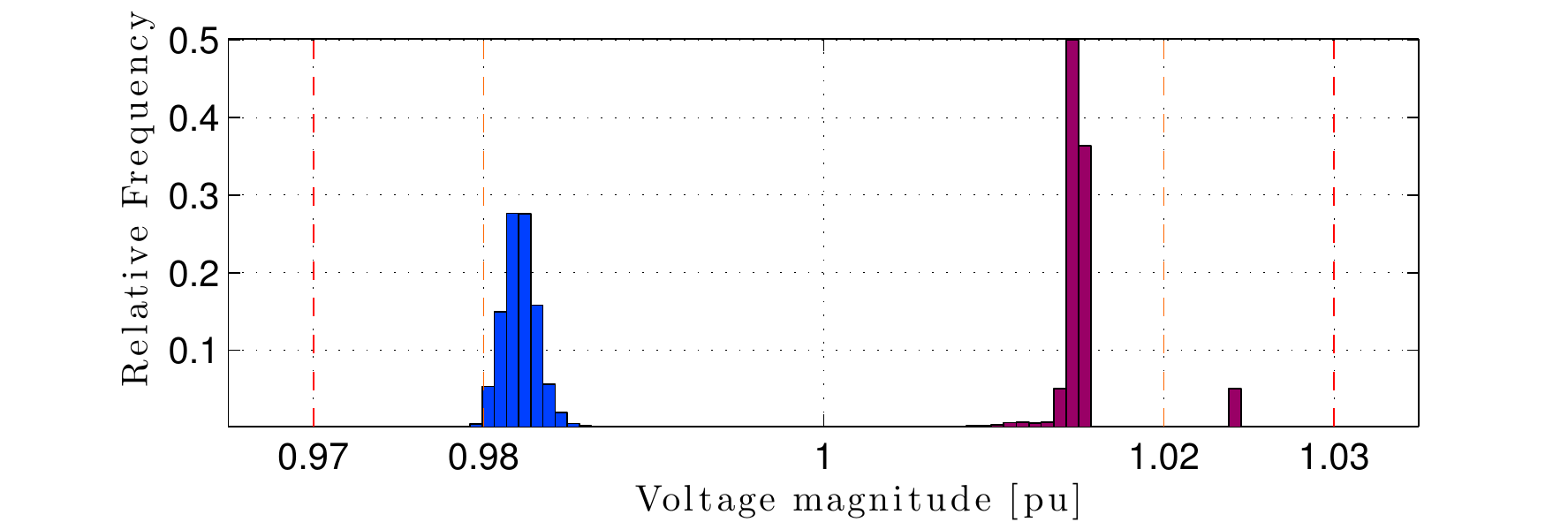}
\vspace*{-1em}
\caption{Histograms of voltage magnitudes on buses 15 and 40 under ADA (top) and PDA (bottom). Dashed lines show regulation limits $\mathcal{V}_A$ and $\mathcal{V}_B$.}\label{f:histograms}
\vspace*{-1em}
\end{figure}

The effect of the average versus the probabilistic constraint on voltage magnitudes was evaluated next. After slow-timescale variables $\bz$ had converged, fast-timescale variables $\by_t$ were calculated for 6,000 instances of $\bxi_t$ using both ADA and PDA. The histograms of the voltage magnitudes on two representative buses are presented in Fig.~\ref{f:histograms}. Under PDA, the average voltage on bus 15 is slightly higher than the average voltage obtained by ADA. In exchange, the instantaneous value of the voltage on bus 15 stays within $\mathcal{V}_A$ with higher probability. A similar behavior is observed for the overvoltage instances on PV bus 40. 

\begin{figure}[t]
\centering
\includegraphics[width=0.95\columnwidth]{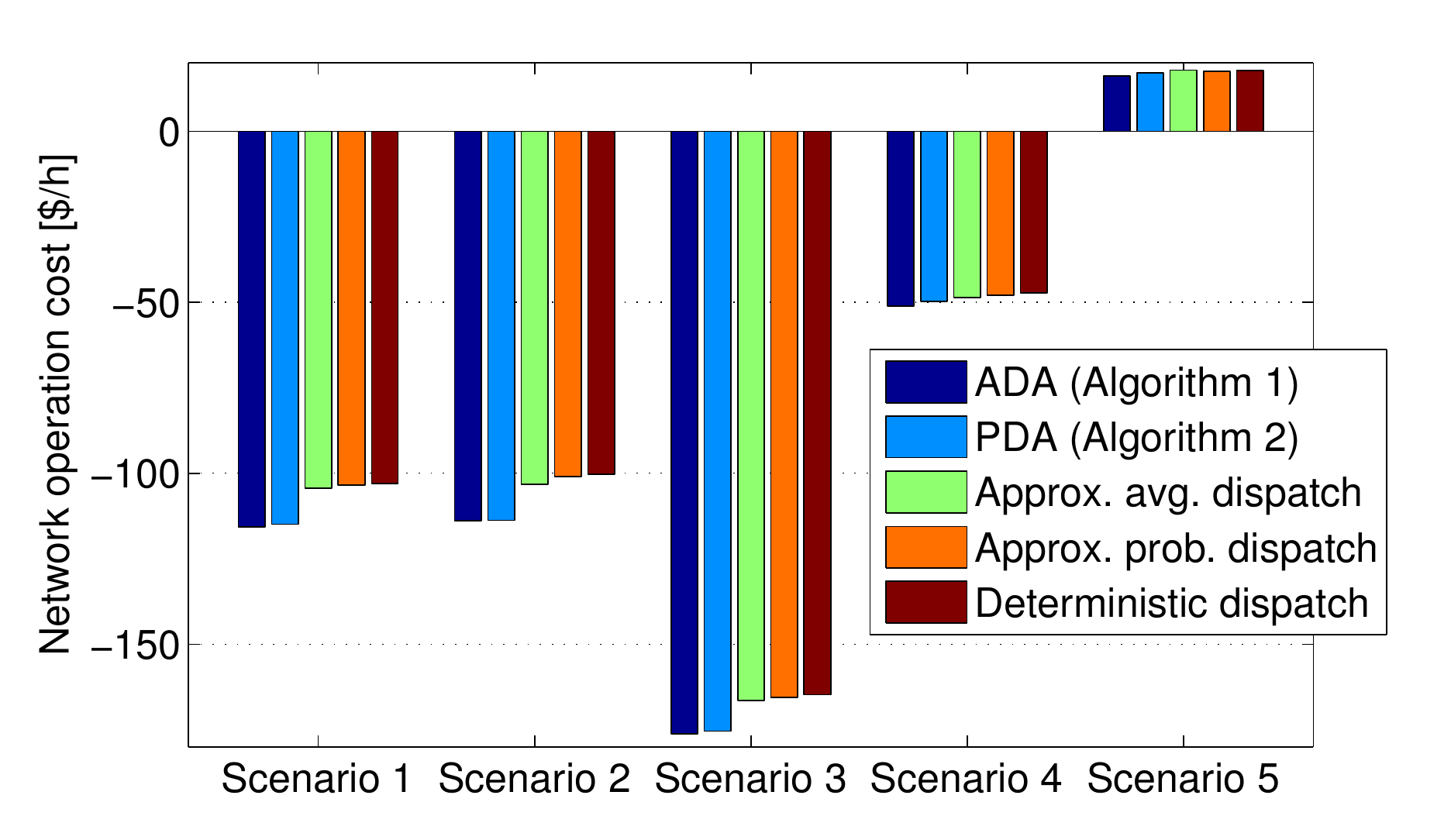}
\vspace*{-0.5em}
\caption{Performance for ADA, PDA, approximate average, approximate probabilistic, and deterministic scheme}.
\vspace*{-1em}
\label{f:barplot}
\end{figure}

ADA and PDA were finally compared to three alternative schemes. The first two, henceforth called \emph{approximate average} and \emph{approximate probabilistic} dispatches, obtained $\bz$ by setting loads and solar generation to their expected values, while variables $\bnu$ were calculated via dual stochastic subgradient, and $\{\by_t\}_{t=1}^T$ were found by solving either \eqref{eqs:inst_opt1} or \eqref{eqs:inst_opt2}, depending on whether the setting is average or probabilistic. The third \emph{deterministic} dispatch found $\bz$ as the approximate schemes do, and $\{\by_t\}_{t=1}^T$ by enforcing $\bv_t \in \mathcal{V}_A$ at all times. Note that the three proposed alternatives provide \textit{feasible} solutions satisfying voltage regulation constraints. The five dispatches were tested under five scenarios: Scenario 1 is the setup described earlier. Scenario 2 involved the tighter voltage limits $\underline{v}_{A}=0.99^2$ and $\overline{v}_{A}=1.01^2$. Scenarios 3, 4, and 5 were generated by scaling the mean value and the standard deviation for loads of scenario 1 by 0.5, 1.5, and 2, respectively. Figure~\ref{f:barplot} shows the expected operation costs for all five scenarios. ADA (PDA) yielded the lowest cost under all scenarios in the average (probabilistic) setting as expected. In all test cases, ADA yielded a slightly lower objective than PDA for $\alpha=0.05$. The loss of optimality entailed by the approximate average and probabilistic schemes is due to the suboptimal choice of $\bz$. The deterministic scheme entailed an additional loss of optimality by preventing the occasional violation of $\mathcal{V}_A$.

\section{Conclusions}\label{s:conclusions}
By nature of renewable generation, electromechanical component limits, and the manner markets operate, energy management of smart distribution grids involves decisions at slower and faster timescales. Since slow-timescale controls remain fixed over multiple PV operation slots, decisions are coupled across time in a stochastic manner. To accommodate solar energy fluctuations, voltages have been allowed to be sporadically overloaded; hence introducing coupling of fast-timescale variables on the average or in probability. Average voltage constraints have resulted in a stochastic convex-concave problem, whereas non-convex probabilistic constraints were tackled using dual decomposition and convex optimization. Efficient algorithms for finding both slow and fast controls using only random samples have been put forth. Our two novel solvers converge in terms of the primal and dual variables, and have attained lower operational costs compared to deterministic alternatives. Although probabilistic constraints have been applied grid-wise, voltages on individual buses remained within limits. Enforcing probabilistic constraints on a per-bus basis, developing decentralized implementations, and including voltage regulators are interesting research directions.

\bibliographystyle{IEEEtran}
\bibliography{myabrv,superpower} 
\end{document}